\documentclass{amsart}

\usepackage{amsthm}
\usepackage{amsmath}
\usepackage{amssymb}
\usepackage{amsfonts}
\usepackage{latexsym}
\usepackage{enumitem}
\usepackage{mathrsfs}
\usepackage[all]{xy} \SelectTips{eu}{}
\usepackage{hyperref}

\newcommand{\numberseries}{\bfseries}   
\newlength{\thmtopspace}                
\newlength{\thmbotspace}                
\newlength{\thmheadspace}               
\newlength{\thmindent}                  
\newtheoremstyle{bfupright head,slanted body}
                {\thmtopspace}{\thmbotspace}
                {\slshape}{\thmindent}{\bfseries}{.}{\thmheadspace}
                {{\numberseries \thmnumber{#2\;}}\thmnote{#3}}

\newtheoremstyle{bfupright head,upright body}
                {\thmtopspace}{\thmbotspace}
                {\upshape}{\thmindent}{\bfseries}{.}{\thmheadspace}
                {{\numberseries \thmnumber{#2\;}}\thmnote{#3}}

\newtheoremstyle{fixed bf head,slanted body}
                {\thmtopspace}{\thmbotspace}{\slshape}
                {\thmindent}{\bfseries}{.}{\thmheadspace}
                {{\numberseries \thmnumber{#2\;}}\thmname{#1}\thmnote{ (#3)}}

\newtheoremstyle{fixed bf head,upright body}
                {\thmtopspace}{\thmbotspace}{\upshape}
                {\thmindent}{\bfseries}{.}{\thmheadspace}
                {{\numberseries \thmnumber{#2\;}}\thmname{#1}\thmnote{ (#3)}}

\newtheoremstyle{numbered paragraph}
                {\thmtopspace}{\thmbotspace}{\upshape}
                {\thmindent}{\upshape}{}{\thmheadspace}
                {{\numberseries \thmnumber{#2.}}}

\theoremstyle{bfupright head,slanted body}
\newtheorem{res}{}[section]             \newtheorem*{res*}{}

\theoremstyle{bfupright head,upright body}
               \newtheorem*{bfhpg*}{}

\theoremstyle{fixed bf head,slanted body}
\newtheorem{thm}[res]{Theorem}          \newtheorem*{thm*}{Theorem}
\newtheorem{prp}[res]{Proposition}      \newtheorem*{prp*}{Proposition}
\newtheorem{cor}[res]{Corollary}        \newtheorem*{cor*}{Corollary}
\newtheorem{lem}[res]{Lemma}            \newtheorem*{lem*}{Lemma}

\theoremstyle{fixed bf head,upright body}
       \newtheorem*{dfn*}{Definition}
           \newtheorem*{rmk*}{Remark}
\newtheorem{exa}[res]{Example}           \newtheorem*{exa*}{Example}

\theoremstyle{numbered paragraph}
\newtheorem{ipg}[res]{}

\newlength{\thmlistleft}        
\newlength{\thmlistright}       
\newlength{\thmlistpartopsep}   
\newlength{\thmlisttopsep}      
\newlength{\thmlistparsep}      
\newlength{\thmlistitemsep}     

\newcounter{eqc}
\newenvironment{eqc}{\begin{list}{\upshape (\textit{\roman{eqc}})}%
    {\usecounter{eqc}%
      \setlength{\leftmargin}{\thmlistleft}%
      \setlength{\labelwidth}{\thmlistleft}%
      \setlength{\rightmargin}{\thmlistright}%
      \setlength{\partopsep}{\thmlistpartopsep}%
      \setlength{\topsep}{\thmlisttopsep}%
      \setlength{\parsep}{\thmlistparsep}%
      \setlength{\itemsep}{\thmlistitemsep}}}%
  {\end{list}}%

\newcounter{prt}
\newenvironment{prt}{\begin{list}{\upshape (\alph{prt})}%
    {\usecounter{prt}%
      \setlength{\leftmargin}{\thmlistleft}%
      \setlength{\labelwidth}{\thmlistleft}%
      \setlength{\rightmargin}{\thmlistright}%
      \setlength{\partopsep}{\thmlistpartopsep}%
      \setlength{\topsep}{\thmlisttopsep}%
      \setlength{\parsep}{\thmlistparsep}%
      \setlength{\itemsep}{\thmlistitemsep}}}%
  {\end{list}}%

\newcounter{rqm}
  {\end{list}}%

\newenvironment{prf*}[1][Proof]{%
  \begin{proof}[\bf #1]
    \setcounter{equation}{0}
    }
  {\end{proof}
}

\newcommand{\pgref}[1]{\ref{#1}}

\renewcommand{\eqref}[1]{(\pgref{eq:#1})}

\newcommand{\thmcite}[2][?]{\cite[thm.~#1]{#2}}
\newcommand{\prpcite}[2][?]{\cite[prop.~#1]{#2}}
\newcommand{\corcite}[2][?]{\cite[cor.~#1]{#2}}
\newcommand{\lemcite}[2][?]{\cite[lem.~#1]{#2}}

\newcommand{\dfncite}[2][?]{\cite[def.~#1]{#2}}

\newcommand{\exacite}[2][?]{\cite[exa.~#1]{#2}}
\newcommand{\eqclbl}[1]{{\upshape(\textit{#1})}}
\newcommand{\proofofimp}[3][:]{\mbox{\eqclbl{#2}$\!\implies\!$\eqclbl{#3}#1}}
\newcommand{\Rop}{R^\circ}

\numberwithin{equation}{res}
\def\urltilda{\kern -.15em\lower .7ex\hbox{\~{}}\kern .04em}

\newcommand{\GPac}{\sf{GPac}}
\newcommand{\GIac}{\sf{GIac}}
\newcommand{\GInj}{\sf{GI}}
\newcommand{\GFac}{\sf{GFac}}
\newcommand{\GPacpd}[2][R]{\operatorname{\GPac-\rm pd}_{#1}#2}
\newcommand{\GIacid}[2][R]{\operatorname{\GIac-\rm id}_{#1}#2}
\newcommand{\GFacpd}[2][R]{\operatorname{\GFac-\rm pd}_{#1}#2}
\newcommand{\Cot}[2][R]{\operatorname{\Cot-\rm id}_{#1}#2}
\newcommand{\GPacgd}[1]{\operatorname{\GPac-\rm gldim}#1}
\newcommand{\GIacgd}[1]{\operatorname{\GIac-\rm gldim}#1}
\newcommand{\GFacgd}[1]{\operatorname{\GFac-\rm gldim}#1}

\newcommand{\Ggd}[1]{\operatorname{Ggldim}#1}
\newcommand{\catD}[1]{\mathsf{D}^\textnormal{b}(#1)}
\newcommand{\catK}[1]{\mathsf{K}^\textnormal{b}(#1)}
\newcommand{\M}[1]{\mathsf{Mod}(#1)}

\newcommand{\Ch}{\mathsf{Ch}}
\newcommand{\lra}{\longrightarrow}
\newcommand{\xra}[2][]{\xrightarrow[#1]{\:#2\:}}
\newcommand{\dif}[2][]{{\partial}^{#2}_{#1}}
\newcommand{\ZZ}{\mathbb{Z}}
\renewcommand{\H}[2][]{\operatorname{H}_{#1}(#2)}

\newcommand{\Z}{\mbox{\rm Z}}

\newcommand{\D}{\mbox{\rm D}}
\newcommand{\Gacpd}{\mbox{\rm Gac-pd}}
\newcommand{\Gacid}{\mbox{\rm Gac-id}}
\newcommand{\Hom}{\mbox{\rm Hom}}
\newcommand{\Ext}{\mbox{\rm Ext}}
\newcommand{\Tor}{\mbox{\rm Tor}}
\newcommand{\id}{\mbox{\rm id}}
\newcommand{\pd}{\mbox{\rm pd}}
\newcommand{\fd}{\mbox{\rm fd}}

\newcommand{\Gpd}{\mbox{\rm Gpd}}
\newcommand{\Gid}{\mbox{\rm Gid}}

\newcommand{\Gacfd}{\mbox{\rm Gac-fd}}
\newcommand{\Ho}{\mbox{\rm Ho}}
\newcommand{\im}{\mbox{\rm Im}}
\newcommand{\Ker}{\mbox{\rm Ker}}
\newcommand{\Coker}{\mbox{\rm Coker}}
\newcommand{\silac}{\mbox{\rm silac}}
\newcommand{\spll}{\mbox{\rm spll}}
\newcommand{\spli}{\mbox{\rm spli}}
\newcommand{\silp}{\mbox{\rm silp}}
\newcommand{\Thb}[2]{#2_{{\scriptscriptstyle\ge}#1}}

   \def\soft#1{\leavevmode\setbox0=\hbox{h}\dimen7=\ht0\advance
    \dimen7 by-1ex\relax\if t#1\relax\rlap{\raise.6\dimen7
    \hbox{\kern.3ex\char'47}}#1\relax\else\if T#1\relax
    \rlap{\raise.5\dimen7\hbox{\kern1.3ex\char'47}}#1\relax
    \else\if d#1\relax\rlap{\raise.5\dimen7\hbox{\kern.9ex
    \char'47}}#1\relax\else\if D#1\relax\rlap{\raise.5\dimen7
    \hbox{\kern1.4ex\char'47}}#1\relax\else\if l#1\relax
    \rlap{\raise.5\dimen7\hbox{\kern.4ex\char'47}}#1\relax
    \else\if L#1\relax\rlap{\raise.5\dimen7\hbox{\kern.7ex
    \char'47}}#1\relax\else\message{accent \string\soft
    \space #1 not defined!}#1\relax\fi\fi\fi\fi\fi\fi}

\begin{document}

\title[Relative global dimensions and stable homotopy categories]%
{Relative global dimensions and stable homotopy categories}

\author[L. Liang]{Li Liang}

\address{L.L. School of Mathematics and Physics, Lanzhou Jiaotong University, Lanzhou 730070, China}

\email{lliangnju@gmail.com}

\urladdr{https://sites.google.com/site/lliangnju}

\author[J. Wang]{Junpeng Wang}

\address{J.W. College of Economics, Northwest Normal University, Lanzhou 730070, China}

\email{wangjunpeng1218@163.com}


\thanks{L. Liang was partly supported by the National Natural Science Foundation of China (Grant No. 11761045), the Foundation of A Hundred Youth Talents Training Program of Lanzhou Jiaotong University, and the Natural Science Foundation of Gansu Province (Grant No. 18JR3RA113); J. Wang was partly supported by the Innovation Foundation of Higher Education of Gansu Province (Grant No. 2020B-087) and a grant from Northwest Normal University (Grant No. NWNU-LKQN2019-14).}




\begin{abstract}
In this paper we study the finiteness of global Gorenstein AC-homological dimensions for rings, and answer the questions posed by Becerril, Mendoza, P\'{e}rez and Santiago. As an application, we show that any left (or right) coherent and left Gorenstein ring has a projective and injective stable homotopy category, which improves the known result by Beligiannis.
\end{abstract}

\maketitle

\thispagestyle{empty}
\section{Introduction}
\noindent
Throughout this work, all rings are assumed to be associative. Let $R$ be a ring; we adopt the convention that an $R$-module is a left $R$-module, and we refer to right $R$-modules as modules over the opposite ring $\Rop$.

Building from Auslander and Bridger's work \cite{MAsMBr69} on modules of finite G-dimension, Enochs, Jenda and Torrecillas \cite{EEnOJn95b,EJT-93} introduced and studied Gorenstein projective, Gorenstein injective and Gorenstein flat modules, and developed ``Gorenstein homological algebra". Such a relative homological algebra theory has been developed rapidly
during the past several years and becomes a rich theory; we refer the reader to, for example, \cite{DB09,BM10,IEm12,EEnOJn95b,EJT-93,Gil17,Wang17} for related works.

For a quasi-Frobenius ring $R$, the category ${\sf Mod}$ of $R$-modules is
a Frobenius category with projective-injective objects all projective (or injective) $R$-modules. So the stable category $\underline{{\sf Mod}}$ modulo projectives is a triangulated category. Furthermore, it is compactly generated; see Krause \cite[sec. 1.5]{Krause00}. It is well known that over an arbitrary ring $R$ the subcategory ${\sf GP}$ (resp., ${\sf GI}$) of Gorenstein projective (resp., Gorenstein injective) $R$-modules is a Frobenius category with projective-injective objects projective (resp. injective) $R$-modules. Also, from a theorem by Christensen, Estrada and Thompson \thmcite[4.5]{CET19} the subcategory ${\sf GF}\cap{\sf Cot}$ of Gorenstein flat and cotorsion $R$-modules is a Frobenius category with projective-injective objects all flat and cotorsion $R$-modules. Hence, the stable categories $\underline{\sf GP}$, $\underline{\sf GI}$ and $\underline{{\sf GF}\cap{\sf Cot}}$ are triangulated categories. It is a natural question when these stable categories are compactly generated. It follows from Beligiannis \cite[lem. 6.6 and thm. 6.7]{Be01} that if $R$ is a right coherent and left perfect or left Morita ring with $\Ggd R<\infty$ then $\underline{\sf GP}\simeq\underline{\sf GI}$ are compactly generated. The same conclusion holds if $R$ is Iwanaga-Gorenstein; see Hovey \thmcite[9.4]{Ho02} or Chen \thmcite[4.1]{Chen11}. One of the main results in this paper is the next improved result; see Corollaries \ref{main-pd and fd-coherent 1} and \ref{main-pd and fd-coherent 2}.

\begin{thm}\label{thm2}
Let $R$ be a ring with $\Ggd R<\infty$.
\begin{prt}
\item If $R$ is right coherent, then $\underline{\sf GP}\simeq\underline{\sf GI}\simeq \underline{\sf GF\cap\sf{Cot}}$ are compactly generated.
\item If $R$ is left coherent, then $\underline{\sf GP}\simeq\underline{\sf GI}$ are compactly generated.
\end{prt}
\end{thm}
\noindent
Here $\Ggd R$ is the global Gorenstein dimension, which is defined as $\Ggd R=\sup\{\Gpd_R M\ |\ M\ \text{is an}\ R\text{-module}\}$. We notice that $R$ satisfies $\Ggd R<\infty$ if and only if $R$ is left Gorenstein\footnote{From Beligiannis \cite{Beligiannis2000}, a ring $R$ is called left Gorenstein if any projective $R$-module has finite injective dimension and any injective $R$-module has finite projective dimension.} ; see \ref{spli}. So as an immediate consequence of Theorem \ref{thm2}, by using \lemcite[6.6]{Be01}, we get the next result that improves \thmcite[6.7]{Be01} by removing the assumption that the ring $R$ should be left perfect or left Morita.

\begin{cor}\label{cor2}
Any left (or right) coherent and left Gorenstein ring has a projective and injective stable homotopy category.
\end{cor}
\noindent
We refer the reader to \cite[def. 6.2]{Be01} for the definition of projective/injective stable homotopy category.

An example is given to show that coherent rings of finite global Gorenstein dimension (or equivalently, left Gorenstein) may not be Iwanaga-Gorenstein nor perfect nor Morita; see Example \ref{exa-2}.

According to \thmcite[6.9]{Beligiannis2000} and \thmcite[4.1]{IEm12}, the finiteness of $\Ggd R$ can be characterized by the existence of the triangulated equivalences $\underline{\sf GP}\simeq \catD{R}/\catK{\sf Prj}$ and/or $\underline{\sf GI}\simeq \catD{R}/\catK{\sf Inj}$, where $\catD{R}$ denotes the bounded derived category of $R$, and $\catK{\sf Prj}$ (resp., $\catK{\sf Inj}$) denotes the bounded homotopy category of projective (resp., injective) $R$-modules. The Verdier quotient triangulated category $\catD{R}/\catK{\sf Prj}$ was first studied by Buchweitz \cite{ROB86} under the name of ¡°stable derived category¡±; it is named by ``singularity category" to emphasize certain homological singularity of the ring $R$ reflected by this quotient category (see Orlov \cite{Orlov04} and Chen \cite{Chen11}). As another immediate consequence of Theorem \ref{thm2}, we see that the singularity categories $\catD{R}/\catK{\sf Prj}\simeq\catD{R}/\catK{\sf Inj}$ are compactly generated over left (or right) coherent rings of finite global Gorenstein dimension (or equivalently, left Gorenstein); see Corollaries \ref{main-pd and fd-coherent 1} and \ref{main-pd and fd-coherent 2}.
\begin{equation*}
  \ast \ \ \ast \ \ \ast
\end{equation*}
\noindent
We prove Theorem \ref{thm2} above by using the finiteness of global Gorenstein AC-homological dimensions.

Gorenstein AC-projective (resp., Gorenstein AC-injective) dimension is defined in terms of resolutions by Gorenstein AC-projective (resp., Gorenstein AC-injective) modules that were initially introduced by Bravo, Gillespie and Hovey \cite{BGH14} as a natural way to extend the notion of Gorenstein projective (resp., Gorenstein injective) modules. We let $\GPacgd R$ and $\GIacgd R$ denote the global Gorenstein AC-projective and global Gorenstein AC-injective dimension of $R$, respectively. That is, $\GPacgd R=\{\Gacpd_R M\ |\ M\ \text{is an}\ R\text{-module}\}$ and $\GIacgd R=\{\Gacid_R M\ |\ M\ \text{is an}\ R\text{-module}\}$. Recently, Becerril, Mendoza, P\'{e}rez and Santiago \cite[6.15]{BMPS19} asked under which conditions on $R$ the following statements are true:
\begin{prt}
\item[$\bullet$] All $R$-modules have finite Gorenstein AC-projective dimension.
\item[$\bullet$] Any $R$-module has finite Gorenstein AC-projective dimension if and only if it has finite Gorenstein AC-injective dimension.
\end{prt}

In Section \ref{Gacgd} we focus on the above two questions. Our main results in this section are the next two theorems, where the first one is used in the proof of Theorem \ref{thm2}.

\begin{thm}
Let $R$ be a ring with $\Ggd R<\infty$.
\begin{prt}
\item If $R$ is right coherent, then $\GPacgd R<\infty$.
\item If $R$ is left coherent, then $\GIacgd R<\infty$.
\end{prt}
\end{thm}
\noindent
This result is proved in Theorem \ref{exa-pd}. The converses of the above statements are not true in general; see Example \ref{exa-1}.

\begin{thm}
If $R$ is a commutative ring, then $\GIacgd R=\GPacgd R$.
\end{thm}
\noindent
This result is proved in Corollary \ref{gacgdim}. However, to the best of our knowledge, we don't know whether the equality $\GIacgd R=\GPacgd R$ holds for an arbitrary ring $R$.

\section{Preliminaries}\label{pre}
\noindent
We begin with some notation and terminology for use throughout this paper.

\begin{ipg}
By an $R$-complex $M$ we mean a complex of $R$-modules as follows:
  \begin{equation*}
    \cdots \lra M_{i+1} \xra{\dif[i+1]{M}} M_i \xra{\dif[i]{M}}
    M_{i-1} \lra \cdots\: .
  \end{equation*}
We frequently (and without warning) identify $R$-modules with $R$-complexes concentrated in degree $0$. For an $R$-complex $M$, we set $\sup M=\sup\{i\in\ZZ\ |\ M_{i}\neq 0\}$ and $\inf M=\inf\{i\in\ZZ\ |\ M_{i}\neq 0\}$. An $R$-complex $M$ is called \emph{bounded} if $\sup M<\infty$ and $\inf M > -\infty$. The symbol $\H[n]{M}$ denotes the $n$th \emph{homology} of $M$, i.e., $\Ker\dif[n]{M}/\im\dif[n+1]{M}$. An $R$-complex $M$ is called \emph{homology bounded} if $\sup\H{M}<\infty$ and $\inf\H{M}> -\infty$. For an $R$-complex $M$, the symbol $M_{\leq n}$ denotes the subcomplex of $M$ with $(M_{\leq n})_i=M_i$ for $i\leq n$ and $(M_{\leq n})_i=0$ for $i > n$, and the symbol $M_{\ge n}$ denotes the quotient complex of $M$ with $(M_{\ge n})_i=M_i$ for $i\ge n$ and $(M_{\ge n})_i=0$ for $i < n$.

We denote by $\catD{R}$ the bounded derived category of $R$-modules, by $\sf{Prj}$ (resp., $\sf{Inj}$, $\sf{Flat}$, and $\sf{Cot}$) the subcategory of projective (resp., injective, flat, and cotorsion) $R$-modules, and by $\catK{\sf Prj}$ (resp., $\catK{\sf Inj}$, and $\catK{\sf FlatCot}$) the bounded homotopy category of projective (resp., injective, and flat and cotorsion) $R$-modules.
\end{ipg}

\begin{ipg}
An $R$-module $M$ is called \emph{Gorenstein projective} \cite{EEnOJn95b} if there exists an exact sequence $\cdots\to P_1\to P_0\to P_{-1}\to\cdots$ of projective $R$-modules such that $M \cong \Coker(P_{1}\to P_{0})$, and it remains exact after applying the functor $\Hom_{R}(-,P)$ for each projective $R$-module $P$. Dually, one has the definition of Gorenstein injective $R$-modules. An $R$-module $M$ is called \emph{Gorenstein flat} \cite{EJT-93} if there exists an exact sequence $\cdots\to F_1\to F_0\to F_{-1}\to \cdots$ of flat $R$-modules such that
$M\cong\Coker(F_1\to F_0)$, and it remains exact after applying the functor $I\otimes_R-$ for each injective $\Rop$-module $I$. We let $\sf GP$ (resp., $\sf GI$, and $\sf GF$) denote the subcategory of Gorenstein projective (resp., Gorenstein injective, and Gorenstein flat) $R$-modules.

The Gorenstein projective dimension of an $R$-module $M$, $\Gpd_R M$, is defined by declaring that $\Gpd_R M\leq n$ if and only if $M$ has a Gorenstein projective resolution of length $n$, that is, there is an exact sequence $0 \to G_n \to \cdots \to G_0 \to M \to 0$ with each $G_i$ Gorenstein projective. The definition of Gorenstein injective dimension, $\Gid_RM$, can be defined dually. We let $\Ggd R$ denote the global Gorenstein projective dimension of $R$, that is, $\Ggd R=\sup\{\Gpd_R M\ |\ M\ \text{is an}\ R\text{-module}\}$. The next result is proved by Bennis and Mahdou \thmcite[1.1]{BM10}, which is used frequently in the paper.
\end{ipg}

\begin{lem}\label{BMth1.1}
For any ring $R$, $\Ggd R=\sup\{\Gid_R M\ |\ M\ \text{is an}\ R\text{-module}\}$.
\end{lem}

\begin{ipg}\label{spli}
Let $\silp R$ denote the supremum of the injective lengths of projective $R$-modules, and $\spli R$ the supremum of the projective lengths of injective $R$-modules. Since an arbitrary direct sum of projective $R$-modules is projective, the invariant $\silp R$ is finite if and only if every projective $R$-module has finite injective dimension. It follows from Beligiannis and Reiten \thmcite[VII. 2.2]{BR07} that every injective $R$-module has finite projective dimension and $\silp R$ is finite if and only if both $\spli R$ and $\silp R$ are finite. Thus by Emmanouil \thmcite[4.1]{IEm12} one gets that the global Gorenstein projective dimension $\Ggd R$ is finite if and only if $R$ is left Gorenstein.
\end{ipg}

\begin{ipg}
Recall from \cite{BGH14} that an $R$-module $F$ is \emph{type FP$_\infty$} if $F$ has a degree-wise finitely generated projective resolution. An $R$-module $A$ is called \emph{absolutely clean} if $\Ext_{R}^{1}(F,A)=0$ for all $R$-modules $F$ of type FP$_\infty$, and an $\Rop$-module $L$ is called \emph{level} if $\Tor_1^{R}(L,F)=0$ for all $R$-modules $F$ of type FP$_\infty$.

Recall from \cite{BGH14} that an $R$-module $M$ is \emph{Gorenstein AC-projective} if there exists an exact sequence $\cdots\to P_1\to P_0\to P_{-1}\to\cdots$ of projective $R$-modules such that $M \cong \Coker(P_{1}\to P_{0})$, and it remains exact after applying the functor $\Hom_{R}(-,L)$ for each level $R$-module $L$.

Dually, an $R$-module $M$ is called \emph{Gorenstein AC-injective} if there exists an exact sequence $\cdots\to I_1\to I_0\to I_{-1}\to\cdots$ of injective $R$-modules such that $M \cong \Coker(I_{1}\to I_{0})$, and it remains exact after applying the functor $\Hom_R(A,-)$ for each absolutely clean $R$-module $A$.

Recall from \cite{BEI18} that an $R$-module $M$ is \emph{Gorenstein AC-flat} if there exists an exact sequence $\cdots\to F_1\to F_0\to F_{-1}\to \cdots$ of flat $R$-modules such that
$M\cong\Coker(F_1\to F_0)$, and it remains exact after applying the functor $A\otimes_R-$ for each absolutely clean $\Rop$-module $A$.

The symbol $\GPac$ (resp., $\GIac$, and $\GFac$) denotes the subcategory of Gorenstein AC-projective (resp., Gorenstein AC-injective, and Gorenstein AC-flat) $R$-modules. It is easy to see that $\GPac\subseteq{\sf GP}$ and $\GIac\subseteq{\sf GI}$.
\end{ipg}

From \thmcite[A.6]{BGH14}, one has the next lemma.

\begin{lem}\label{gacpf}
All Gorenstein AC-projective $R$-modules are Gorenstein AC-flat. That is, $\GPac\subseteq\GFac$.
\end{lem}

\begin{ipg}\label{global}
The Gorenstein AC-projective dimension of $R$-module $M$, $\Gacpd_R M$, is defined by declaring that $\Gacpd_R M\leq n$ if and only if $M$ has a Gorenstein AC-projective resolution of length $n$, that is, there is an exact sequence $0 \to G_n \to \cdots \to G_0 \to M \to 0$ with each $G_i$ Gorenstein AC-projective. The Gorenstein AC-injective and Gorenstein AC-flat dimensions are defined similarly, which are denoted $\Gacid_R M$ and $\Gacfd_R M$, respectively.

Let $\GPacgd R$ (resp., $\GIacgd R$, and $\GFacgd R$) denote the global Gorenstein AC-projective (resp., global Gorenstein AC-injective, and global Gorenstein AC-flat) dimension of $R$. For example,
$$\GPacgd R=\{\Gacpd_R M\ |\ M\ \text{is an}\ R\text{-module}\}.$$
By Lemma \ref{BMth1.1}, one has
\begin{equation}\label{2.7.1}
\Ggd R\leq\min\{\GPacgd R, \GIacgd R\}\:.
\end{equation}
\end{ipg}

\begin{ipg}
A pair $(\sf{X},\sf{Y})$ of subcategories of $R$-modules is called a \emph{cotorsion pair} if
$\sf{X}^\perp=\sf{Y}$ and $\sf{Y}={^\perp\sf{X}}$. Here ${\sf X}^\perp=\{ A~|~\Ext^1_{R}(X,A)=0\ \text{for all}\ X\in\sf{X}\}$, and similarly one can define ${^\perp\sf{X}}$.
A cotorsion pair $(\sf{X},\sf{Y})$ is said to be \emph{hereditary} if $\Ext^n_R(X,Y)=0$ for all $X\in \sf{X}$, $Y\in \sf{Y}$ and $n\geq 1$, or equivalently, if ${\sf Y}$ is injectively coresolving (that is, whenever $0 \to L' \to L \to L'' \to 0$ is exact with $L', L \in\sf{Y}$ then $L''$ is also in $\sf Y$). A cotorsion pair $(\sf{X},\sf{Y})$ is called \emph{complete} if for any $R$-module $A$, there exist exact
sequences $0\to Y\to X\to A\to 0$ and/or $0\to A\to Y'\to X'\to 0$ with $X,X'\in\sf{X}$ and $Y,Y'\in\sf{Y}$.
\end{ipg}

\section{Global Gorenstein AC-homological dimensions}\label{Gacgd}
\noindent
In this section we focus on the global Gorenstein AC-projective/injective dimension. We let $\silac R=\sup\{\id_RM~|~M\ \text{is an absolutely clean}\ R\text{-module}\}$.

\begin{lem}\label{m-abclean}
Let $R$ be a ring. Then there exists an equality
\begin{equation*}
\max\{\Ggd R, \silac R\}=\GIacgd R.
\end{equation*}
\end{lem}
\begin{prf*}
For the inequality ``$\geq$", we let $\max\{\Ggd R, \silac R\}=m<\infty$. Then all absolutely clean $R$-module have finite injective dimension.
This implies that $\GIac=\GInj$. Thus we have $\GIacgd R=\Ggd R\leq m$.

For the inequality ``$\leq$", we let $\GIacgd R= n<\infty$. It is easy to see that $\Ggd R\leq\GIacgd R=n$. Next we prove that $\silac R\leq n$. Let $A$ be an absolutely clean $R$-module. For each $R$-module $M$, one has $\GIacid M\leq n$. So there is an exact sequence $0 \to M \to G^0 \to \cdots \to G^n \to 0$ with each $G^i$ Gorenstein AC-injective. Thus $\Ext^{n+1}_R(A,M)\cong\Ext^{1}_R(A,G^n)=0$. This yields that $A$ has finite projective dimension at most $n$. So $A$ has finite injective dimension at most $n$ by \corcite[2.7]{BM10}, as $\Ggd R\leq\GIacgd R=n$; see (\ref{2.7.1}). Thus one has $\silac R\leq n$.
\end{prf*}

The next result is immediate by Lemma \ref{m-abclean}.

\begin{lem}\label{abclean}
Let $R$ be a ring with $\GIacgd R$ finite. Then all absolutely clean $R$-modules have finite injective dimension at most $\GIacgd R$. Hence, all Gorenstein injective $R$-modules are Gorenstein AC-injective.
\end{lem}

The following two results are proved dually, where we let
$$\spll R=\sup\{\pd_RM~|~M\ \text{is a level}\ R\text{-module} \}.$$

\begin{lem}\label{m-level}
Let $R$ be a ring. Then there exists an equality
\begin{equation*}
\max\{\Ggd R, \spll R\}=\GPacgd R.
\end{equation*}
\end{lem}

\begin{lem}\label{level}
Let $R$ be a ring with $\GPacgd R$ finite. Then all level $R$-modules have finite projective dimension at most $\GPacgd R$. Hence, all Gorenstein projective $R$-modules are Gorenstein AC-projective.
\end{lem}

\begin{thm}\label{GIacgd}
The following statements hold:
\begin{prt}
\item If\: $\GIacgd R<\infty$, then there is an inequality $\GIacgd R\leq\GPacgd R$.
\item If\: $\GIacgd\Rop<\infty$, then there is an inequality $\GPacgd R\leq\GIacgd R$.
\end{prt}
\end{thm}
\begin{prf*}
(a) By Lemma \ref{abclean} all Gorenstein injective $R$-modules are Gorenstein AC-injective, so one has $\GIacgd R=\Ggd R\leq \GPacgd R$.

(b) We may assume that $\GIacgd R$ is finite. By Lemma \ref{abclean} all absolutely clean $\Rop$-modules have finite injective dimension, as $\GIacgd\Rop<\infty$. Let $L$ be a level $R$-module. Then by \thmcite[2.12]{BGH14} $L^+$ is a absolutely clean $\Rop$-module, and so  $\fd_R L=\id_RL^{+}<\infty$. Thus $\pd_RL<\infty$ by \corcite[2.7]{BM10} as $\Ggd R<\infty$.
Hence, one has $\sf{GP}=\GPac$; see Lemma \ref{level}. Thus $\GPacgd R=\Ggd R\leq\GIacgd R$; see \ref{global}.
\end{prf*}

The next result is proved dually.

\begin{thm}\label{GPacgd}
The following statements hold:
\begin{prt}
\item If\: $\GPacgd R<\infty$, then there is an inequality $\GPacgd R\leq\GIacgd R$.
\item If\: $\GPacgd\Rop<\infty$, then there is an inequality $\GIacgd R\leq\GPacgd R$.
\end{prt}
\end{thm}

The next corollary advertised in the introduction is immediate by Theorems \ref{GIacgd} and \ref{GPacgd}.

\begin{cor}\label{gacgdim}
If $R$ is a commutative ring, then $\GIacgd R=\GPacgd R$.
\end{cor}

\begin{lem}\label{Gac p and f4}
Let $M$ be an $R$-module and $n$ a nonnegative integer. Then the following conditions are equivalent.
\begin{eqc}
 \item $\GFacpd M\leq n$.
 \item There is an exact sequence $0\to M \to F \to N \to 0$ of $R$-modules
with $\fd_R F\leq n$ and $N\in\GPac$.
\end{eqc}
\end{lem}
\begin{prf*}
\proofofimp{i}{ii} We prove the result by induction on $n$. The case where $n=0$ holds by Proposition \ref{Gac p and f}.
Now let $n>0$. Consider an exact sequence
$0\to K \to H \to M \to 0$ of $R$-modules
with $H$ flat. Then one has $\GFacpd K\leq n-1$, and so by induction,
there is an exact sequence
$0\to K\to H'\to G\to 0$ of $R$-modules
with $\fd_RH'\leq n-1$ and $G\in \GPac$.
Consider the following pushout diagram
$$\xymatrix@C=20pt@R=20pt{&0\ar[d] &0\ar[d]\\
  0 \ar[r] & K \ar[r]^{}\ar[d] & H\ar[r]^{}\ar[d] & M \ar[r]\ar@{=}[d] & 0\\
  0 \ar[r] & H' \ar[r]^{}\ar[d] & H'' \ar[r]^{}\ar[d] & M \ar[r] & 0.\\
  &G\ar@{=}[r]\ar[d]&G\ar[d]\\
  &0&0}$$
In the middle columnn, by Lemma \ref{gacpf} both $H$ and $G$ are in $\GFac$, so is $H''$.
Whence, by Proposition \ref{Gac p and f}, there is an exact sequence
$0\to H''\to L\to N\to 0$ of $R$-modules
with $L\in \mathsf{Flat}$ and $N\in \GPac$.
Now we obtain another pushout diagram
$$\xymatrix@C=20pt@R=20pt{&&0\ar[d] &0\ar[d]\\
  0 \ar[r] & H' \ar[r]^{}\ar@{=}[d] & H''\ar[r]^{}\ar[d] & M \ar[r]\ar[d] & 0\\
  0 \ar[r] & H' \ar[r]^{} & L \ar[r]^{}\ar[d] & F \ar[r]\ar[d] & 0.\\
  &&N\ar@{=}[r]\ar[d]&N\ar[d]\\
  &&0&0}$$
Since, in the middle row, $L\in \mathsf{Flat}$ and $\fd_RH'\leq n-1$, it follows that $\fd_RF\leq n$. So the condition \eqclbl{ii} holds by the rightmost non-zero column.

\proofofimp[]{ii}{i} Assume that there is an exact sequence $0\to M \to F \to N \to 0$ of $R$-modules
with $\fd_R F\leq n$ and $N\in\GPac$. Since $(\GPac,{\GPac}^{\perp})$ is a hereditary complete cotorsion pair by Gillespie \cite[fact. 10.2]{Gil16}, there is an exact sequence $0\to E \to L \to F \to 0$ of $R$-modules
with $L\in \GPac$ and $E\in{\GPac}^{\perp}$.
Consider the next pullback diagram
$$\xymatrix@C=20pt@R=20pt{&0\ar[d] &0\ar[d]\\
  & E \ar[d]^{}\ar@{=}[r]& E \ar[d]^{}\\
  0 \ar[r] & Q \ar[r]^{} \ar[d]^{}& L \ar[r]^{}\ar[d] & N \ar[r]\ar@{=}[d] & 0\\
 0 \ar[r] & M \ar[r]^{} \ar[d]^{}& F \ar[r]^{}\ar[d] & N \ar[r] & 0.\\
  &0&0}$$
By the middle column one gets $\GFacpd E\leq n-1$ since $\GFacpd F\leq \fd_RF\leq n$ and $L\in \GPac\subseteq\GFac$; see Lemma \ref{gacpf}.
By the middle row one gets that $Q$ is in $\GPac\subseteq\GFac$ since $N$ and $L$ are in $\GPac$.
Thus, by the first non-zero column one has $\GFacpd M\leq n$.
\end{prf*}

\begin{prp}\label{Gfacgf}
Let $R$ be a ring with $\GFacgd R<\infty$. Then all Gorenstein flat $R$-modules are Gorenstein AC-flat.
\end{prp}
\begin{prf*}
We assume that $\GFacgd R=n<\infty$. Let $M$ be a Gorenstein flat $R$-module. Then one has $\GFacpd M\leq n$, and so by Lemma \ref{Gac p and f4} there is an exact sequence $0\to M \to F \to N \to 0$ of $R$-modules with $\fd_R F\leq n$ and $N\in\GPac\subseteq\GFac$. Since $M$ and $N$ are Gorenstein flat, a recent result by \v{S}aroch and \v{S}\v{t}ov\'{\i}\v{c}ek \thmcite[3.11]{SS18} yields that $F$ is Gorenstein flat. Thus $F$ is flat, and hence $M$ is Gorenstein AC-flat by Proposition \ref{Gac p and f}.
\end{prf*}

\begin{ipg}\label{cot-id}
The cotorsion dimension of $R$-module $M$, ${\sf Cot}{\text-}\id_R M$, is defined by declaring that ${\sf Cot}{\text-}\id_R M\leq n$ if and only if $M$ has a cotorsion coresolution of length $n$, that is, there is an exact sequence $0 \to M \to C^0 \cdots \to C^{n} \to 0$ with each $C^i$ cotorsion. We let $\mathsf{Cot}{\text-}{\rm gldim}R=\sup\{{\sf Cot}{\text-}\id_R M\ |\ M\ \text{is an}\ R\text{-module}\}$.
\end{ipg}

The next result was proved by Mao and Ding in \thmcite[19.2.14]{MD06}.

\begin{lem}\label{fpd}
For each $R$-module $M$ there exists an inequality
\begin{equation*}
\pd_RM\leq\fd_RM+\mathsf{Cot}{\text-}{\rm gldim}R.
\end{equation*}
\end{lem}

The next result is used in the proofs of Corollaries \ref{main-pd and fd} and \ref{main-pd and fd-coherent 1}.

\begin{thm}\label{cotgldim}
Let $R$ be a ring. Then there exist inequalities
\begin{equation*}
\max\{\GFacgd R, \mathsf{Cot}{\text-}{\rm gldim}R\}\leq\GPacgd R\leq\GFacgd R+ \mathsf{Cot}{\text-}{\rm gldim}R.
\end{equation*}
In particular, $\GPacgd R$ is finite if and only if $\GFacgd R$ and $\mathsf{Cot}{\text-}{\rm gldim}R$ are finite.
\end{thm}
\begin{prf*}
For the first inequality one let $\GPacgd R=n<\infty$. We notice that all Gorenstein AC-projective modules are Gorenstein AC-flat. So one has $\GFacgd R$ $\leq n$. Let $F$ be an flat $R$-module. Since $\Ggd R\leq\GPacgd R=n$, one has $\pd_RF\leq n$ by \corcite[2.7]{BM10}. Thus \corcite[7.2.6]{MD06} yields $\mathsf{Cot}{\text-}{\rm gldim}R\leq n$.

For the second inequality we let $\GFacgd R=n<\infty$ and $\mathsf{Cot}{\text-}{\rm gldim}R=m<\infty$. Let $M$ be an $R$-module. Then $\GFacpd M\leq n$. By Lemma \ref{Gac p and f4}, one gets an exact sequence $0\to M \to F \to N \to 0$ of $R$-modules
with $\fd_RF\leq n$ and $N \in\GPac$. So $\pd_RH\leq n+m$; see Lemma \ref{fpd}. Similar to the proof of \proofofimp[]{ii}{i} in Lemma \ref{Gac p and f4} one gets $\GPacpd M\leq n+m$. Thus $\GPacgd R\leq n+m$.
\end{prf*}

Next we give some rings that have finite global Gorenstein AC-projective/injective dimension.

\begin{thm}\label{exa-pd}
For a ring $R$ with $\Ggd R<\infty$, the following statements hold:
\begin{prt}
\item If $R$ is left coherent, then $\GIacgd R<\infty$.
\item If $R$ is right coherent, then $\GPacgd R<\infty$.
\end{prt}
\end{thm}
\begin{prf*}
(a) By Lemma \ref{m-abclean}, it is suffices to show that $\silac R<\infty$.
Let $A$ be an absolutely clean $R$-module, and let $\Ggd R=n<\infty$.
Then \corcite[2.9]{BGH14} yields that $A$ is FP-injective since $R$ is left coherent.
Hence, there is a pure exact sequence $0\to A\to I\to C\to 0$ of $R$-modules with $I$ injective. By \corcite[2.7]{BM10}, one has $\fd_RI\leq n$. It follows that $\fd_RA\leq \fd_RI\leq n$,
and hence one has $\id_RA\leq n$ again by \corcite[2.7]{BM10}.
This gives that $\silac R\leq n<\infty$.

(b) By Lemma \ref{m-level}, it is suffices to show that $\spll R<\infty$.
Let $L$ be a level $R$-module, and let $\Ggd R=n<\infty$. Then \corcite[2.11]{BGH14} yields that $L$ is flat since $R$ is right coherent.
Hence one has $\pd_RL\leq n$ by \corcite[2.7]{BM10}. This gives that $\spll R\leq n<\infty$.
\end{prf*}

In the followin we give an example to show that the converses of the statements in Theorem \ref{exa-pd} are not true in general. Before that we give some facts.

\begin{ipg}\label{product}
Let $R=\prod_{i=1}^{n}R_i$ be a direct product of rings. If $M_i$ is an $R_i$-module for $i=1,2,\cdots,n$ then $M=M_1 \oplus M_2 \oplus \cdots \oplus M_n$ is an $R$-module. Conversely, if $M$ is an $R$-module then it is of the form $M=M_1 \oplus M_2 \oplus \cdots \oplus M_n$, where $M_i$ is an $R_i$-module for $i=1,2,\cdots,n$. It is easy to see that the following equalities hold
\begin{equation*}
\GPacpd M = \sup\{\GPacpd[R_i]M_i~|~i = 1, ..., n\}
\end{equation*}
and
\begin{equation*}
\GIacid M = \sup\{\GIacid[R_i]M_i~|~i = 1, ..., n\},
\end{equation*}
which are parallel to the well-known ones about projective and injective dimension, respectively.
So one gets that $R$ is of finite global Gorenstein AC-projective/injective dimension if and only if each $R_i$ is so; the same conclusion holds for global dimension.
On the other hand, it is known that $R$ is left/right coherent if and only if each $R_i$ is so.
\end{ipg}

\begin{exa} \label{exa-1}
Let $R=D+(x_1,x_2)K[x_1, x_2]$, where $D$ is a Dedekind domain and $K$ its
quotient field. According to Kirkman and Kuzmanovich \cite[Example in p.128]{KK88},
$R$ is a commutative non-coherent ring of finite global dimension.
On the other hand, there exists a commutative Iwanaga-Gorenstein ring $S$ of infinite global dimension; see Bennis \cite[p.857]{DB09}. So $S$ has finite global Gorenstein AC-projective dimension and finite global Gorenstein AC-injective dimension; see Theorem \ref{exa-pd}.
Hence, $R\times S$ has finite global Gorenstein AC-projective dimension and finite global Gorenstein AC-injective dimension.
However, $R\times S$ is neither of finite global dimension nor coherent.
\end{exa}

\section{Compactly generatedness of singularity categories}
\noindent
We now turn to study the compactly generatedness of singularity categories and stable categories with respect to Gorenstein AC-homological modules, and prove Theorem \ref{thm2} advertised in the introduction. We open this section with the following terminology.

\begin{ipg}\label{dim comp}
Let $(\sf{A}, \sf{B})$ be a cotorsion pair in $\M{R}$, and let $X$ be an $R$-complex. From Yang and Ding \cite{YD15}, the \emph{$\sf{A}$-projective dimension} of $X$, ${\sf A}{\text-}\pd_RX$, is defined as
\begin{equation*}
{\sf A}{\text-}\pd_RX= \inf\{\sup A~|~X \simeq A\ \mathrm{in}\ \D(R)~\text{with}~ A \in \sf{dgA}\}.
\end{equation*}
The \emph{$\sf{B}$-injective dimension} of $X$, ${\sf B}{\text-}\id_RX$, is defined as
\begin{equation*}
{\sf B}{\text-}\id_RX= \inf\{-\inf B~|~X \simeq B\ \mathrm{in}\ \D(R)~\text{with}~ B \in \sf{dgB}\}.
\end{equation*}
Here $\sf dgA$ and $\sf dgB$ denote the subcategories of dg-${\sf A}$ complexes and dg-$\sf B$ complexes, respectively; see Gillespie \cite{Gi04}.

It is known that $({\sf Flat}, {\sf Cot})$ is a complete hereditary cotorsion pair, and by \cite{BGH14} and Proposition \ref{Gac p and f1} $(\GPac, \GPac^{\perp})$, $({^{\perp}\GIac},\GIac)$ and $(\GFac, (\GPac)^{\perp}\cap \sf{Cot})$ are complete hereditary cotorsion pairs. So for an $R$-complex $X$ we have the definitions of ${\GPac}{\text-}\pd_RX$, ${\GFac}{\text-}\pd_RX$, ${\GIac}{\text-}\id_RX$ and ${\sf Cot}{\text-}\id_RX$. We let $\catD{R}_{\widehat{\GPac}}$ (resp., $\catD{R}_{\widehat{\GFac}\cap\widehat{\sf{Cot}}}$, and $\catD{R}_{\widehat{\GIac}}$) denote the the triangulated subcategory of $\catD{R}$
consisting of all homology bounded complexes $X$ with ${\GPac}{\text-}\pd_RX<\infty$ (resp, ${\GFac}{\text-}\pd_RX<\infty$ and ${\sf Cot}{\text-}\id_RX<\infty$, and ${\GIac}{\text-}\id_RX<\infty$). It is easy to see that for an $R$-module $M$ (viewed as an $R$-complex concentrated in degree $0$), the definitions of ${\GPac}{\text-}\pd_RM$, ${\GFac}{\text-}\pd_RM$, ${\GIac}{\text-}\id_RM$ and ${\sf Cot}{\text-}\id_RM$ are the same as in \ref{global} and \ref{cot-id}.
\end{ipg}

\begin{lem}\label{frobe}
The subcategory $\GPac$ $($resp., $\GFac\cap\sf{Cot}$, and $\GIac$$)$
together with all short exact sequences in $\GPac$ $($resp., $\GFac\cap\sf{Cot}$, and $\GIac$$)$
forms a Frobenius category with projective-injective objects all projective $($resp., flat-cotorsion, and injective$)$ $R$-modules.
\end{lem}
\begin{prf*}
We give a straight proof for the case $\GIac$; see \ref{Frobenius exact category} for the other ones.

The subcategory $\GIac$, together with all short exact sequences in $\GIac$, forms an exact category, as $\GIac$ is closed under extensions by \lemcite[5.6]{BGH14}.

For $I\in {\sf Inj}$ and $G\in \GIac$, one gets that $\Ext_R^{1}(G,I)=0=\Ext_R^{1}(I,G)$, which yields that all injective $R$-modules are both projectives and injectives in $\GIac$. Conversely, let $M$ (resp., $N$) be a injective (resp., projective) object in $\GIac$. Then there exist split exact sequences $0\to M\to I\to M' \to 0$ and $0\to N'\to H\to N\to 0$ with $I,H\in {\sf Inj}$ and $M',N'\in \GIac$. So both $M$ and $N$ are in $\sf Inj$. Thus projectives and injectives in $\GPac$ are exactly injective $R$-modules.

Finally, for every $G\in\GIac$ there exist exact sequences $0\to G\to I'\to G'\to0$ and $0\to G''\to I''\to G\to0$ with $I', I''\in {\sf Inj}$ and $G', G''\in \GIac$, so the subcategory $\GIac$ has enough injectives and enough projectives.
\end{prf*}

\begin{ipg}
By Lemma \ref{frobe}, the stable category $\underline{\GPac}$ (resp., $\underline{\GFac\cap\sf{Cot}}$, and $\underline{\GIac}$) modulo projectives (resp., flat-cotorsions, and injectives) is a triangulated category.
\end{ipg}

\begin{thm}\label{equivalence}
The following conditions are equivalent.
 \begin{eqc}
 \item $\GPacgd R<\infty$.
 \item There is an equality $\catD{R}_{\widehat{\GPac}}=\catD{R}$.
 \item The natural functor $\mathrm{F}: \underline{\GPac} \to \catD{R}/\catK{\sf Prj}$ induced by the compositions
     \begin{equation*}
     \sf GPac\hookrightarrow \catD{R}_{\widehat{\GPac}}\to\catD{R}_{\widehat{\GPac}}/\catK{\sf Prj}\hookrightarrow\catD{R}/\catK{\sf Prj}
     \end{equation*}
     is a triangulated equivalence.
\end{eqc}
\end{thm}
\begin{prf*}
\proofofimp{i}{ii} Fix $P\in\catD{R}$. It suffices to show that ${\GPac}{\text-}\pd_RP<\infty$. Without loss of generality, we may assume that $P$ is bounded as follows:
\begin{equation*}
P=0\to P_k \to P_{k-1} \to \cdots \to P_{1}\to P_{0}\to 0\:.
\end{equation*}
Consider the exact sequence $0 \to P_0 \to P \to \Thb{1}{P} \to 0$ of $R$-complexes. Since $P_0$ and $\Thb{1}{P}$ have finite Gorenstein AC-projective dimension by \eqclbl{i} and induction on $k$, respectively, so does $P$.

\proofofimp{ii}{i} Each $R$-module $M$, viewed as an $R$-complex concentrated in degree $0$, is in $\catD{R}$. Thus $M\in\catD{R}_{\widehat{\GPac}}$,
and so $\GPacpd{M}<\infty$; see \ref{dim comp}. Note that $(\GPac, \GPac^{\perp})$ is a complete hereditary cotorsion pair. It is a standard way to see that for any family $(M_i)_{i\in\Lambda}$ of $R$-modules there is an equality
$$\GPacpd{(\oplus_{i\in\Lambda} M_i)}=\sup\{\GPacpd{M_i}\ |\ i\in\Lambda\}.$$
Thus it is easy to verify that the condition \eqclbl{i} holds.

\proofofimp{ii}{iii} From a result by Di, Liu, Yang and Zhang \corcite[5.9]{DLYZ18}, the induced natural functor $\mathrm{F}: \underline{\GPac} \to \catD{R}_{\widehat{\GPac}}/\catK{\sf Prj}$ is a triangulated equivalence, so the statement \eqclbl{iii} follows from \eqclbl{ii}.

\proofofimp{iii}{ii} It is clear that $\catD{R}_{\widehat{\GPac}}\subseteq\catD{R}$. Conversely, we let $X\in\catD{R}$ ($X$ is also an object of $\catD{R}/\catK{\sf Prj}$). By \eqclbl{iii} and \corcite[5.9]{DLYZ18}, the functor
\begin{equation*}
\catD{R}_{\widehat{\GPac}}/\catK{\sf Prj}\hookrightarrow\catD{R}/\catK{\sf Prj}
\end{equation*}
is a triangulated equivalence. We notice that each triangulated equivalence is dense. So $X$ is isomorphic to an $R$-complex in $\catD{R}_{\widehat{\GPac}}/\catK{\sf Prj}$. It follows that $X$ is in $\catD{R}_{\widehat{\GPac}}$.
\end{prf*}

\begin{thm}\label{equivalence-f}
The following conditions are equivalent.
 \begin{eqc}
\item $\GFacgd R<\infty$ and $\mathsf{Cot}{\text-}{\rm gldim}R<\infty$.
\item There is an equality $\catD{R}_{\widehat{\GFac}\cap\widehat{\sf{Cot}}}=\catD{R}$.
\item There is a triangulated equivalence
   \begin{equation*}
   \underline{\GFac\cap\sf{Cot}}\simeq \catD{R}/\catK{\sf FlatCot}.
   \end{equation*}
\item The natural functor $\mathrm{F}: \underline{\GFac\cap\sf{Cot}} \to \catD{R}/\catK{\sf FlatCot}$ induced by the compositions
     $\GFac\cap\sf{Cot}\hookrightarrow \catD{R}_{\widehat{\GFac}\cap\widehat{\sf{Cot}}}\to \catD{R}_{\widehat{\GFac}\cap\widehat{\sf{Cot}}}/\catK{\sf FlatCot}\hookrightarrow\catD{R}/\catK{\sf FlatCot}$
     is a triangulated equivalence.
\end{eqc}
\end{thm}
\begin{prf*}
Analogous to the proof of Theorem \ref{equivalence}, using Corollary \ref{Gac p and f3} instead of \corcite[5.9]{DLYZ18}.
\end{prf*}

\begin{cor}\label{main-pd and fd}
Let $R$ be a ring with $\GPacgd R$ finite. Then
\begin{equation*}
\catD{R}/\catK{\sf Inj}\simeq\catD{R}/\catK{\sf Prj}\simeq\underline{\GPac}\simeq \underline{\GFac\cap\sf{Cot}}\simeq \catD{R}/\catK{\sf FlatCot}
\end{equation*}
are compactly generated.
\end{cor}
\begin{prf*}
The first equivalence in the statement holds by \thmcite[6.9]{Beligiannis2000} since $\Ggd R$ is finite; see \ref{2.7.1}. The second equivalence follows from Theorem \ref{equivalence}, the third one holds by Corollary \ref{triequivalence}, and the last one follows from Theorems \ref{cotgldim} and \ref{equivalence-f}. By a careful
reading of the proof of Gillespie \thmcite[6.2]{Gil19}, one gets that $\underline{\GPac}$ is compactly generated.
\end{prf*}

It is from \cite[lem. 6.6 and thm. 6.7]{Be01} that if $R$ is a right coherent and left perfect or left Morita ring with $\Ggd R<\infty$ then $\underline{\sf GP}\simeq\underline{\sf GI}$ are compactly generated. The same conclusion holds if $R$ is Iwanaga-Gorenstein; see \thmcite[9.4]{Ho02} or \thmcite[4.1]{Chen11}. We have the next improved result.

\begin{cor}\label{main-pd and fd-coherent 1}
Let $R$ be a right coherent ring with $\Ggd R<\infty$. Then
\begin{equation*}
\catD{R}/\catK{\sf Inj}\simeq\catD{R}/\catK{\sf Prj}\simeq\underline{\sf GP}\simeq\underline{\sf GI}\simeq \underline{\sf GF\cap\sf{Cot}}\simeq \catD{R}/\catK{\sf FlatCot}
\end{equation*}
are compactly generated.
\end{cor}
\begin{prf*}
By Theorem \ref{exa-pd} one has $\GPacgd R$ finite. We notice that $\underline{\sf GP}\simeq\underline{\sf GI}$ by \thmcite[6.9]{Beligiannis2000} as $\Ggd R<\infty$. On the other hand, by Lemma \ref{level} one has $\GPac={\sf GP}$, and the equality $\GFac={\sf GF}$ holds by Proposition \ref{Gfacgf} and Theorem \ref{cotgldim}. So the desired result in the statement follows from Corollary \ref{main-pd and fd}.
\end{prf*}

Dual to the proof of Theorem \ref{equivalence}, we have the following result.

\begin{thm}\label{equivalence-dual}
The following conditions are equivalent.
 \begin{eqc}
 \item $\GIacgd R<\infty$.
  \item There is an equality $\catD{R}_{\widehat{\GIac}}=\catD{R}$
 \item The natural functor $\mathrm{F}: \underline{\GIac} \to \catD{R}/\catK{\sf Inj}$ induced by the compositions
     \begin{equation*}
     \sf GIac\hookrightarrow \catD{R}_{\widehat{\GIac}}\to\catD{R}_{\widehat{\GIac}}/\catK{\sf Inj}\hookrightarrow\catD{R}/\catK{\sf Inj}
     \end{equation*}
     is a triangulated equivalence.
\end{eqc}
\end{thm}

From Gillespie \dfncite[5.1]{Gi17}, a complex $I$ of injective $R$-modules is called \emph{AC-injective} if each chain map into $I$ from an acyclic complex with each cycle absolutely clean
is null homotopic.

\begin{prp}\label{kac-inj}
Let $R$ be a ring with $\GIacgd R$ finite. Then all complexes of injective $R$-modules are AC-injective.
\end{prp}
\begin{prf*}
We let $\sf{dw}\widetilde{\sf Inj}$ denote the subcategory of complexes of injective $R$-modules. Let $I\in\sf{dw}\widetilde{\sf Inj}$, and let $\alpha: X\to I$ be an homomorphisms of $R$-complexes with $X$ acyclic and each cycle $\Z_{i}(X)$ absolutely clean. Next we prove that $\alpha$ is null homotopic. Set $n=\GIacgd R<\infty$. Then by Lemma \ref{abclean} each cycle $\Z_{i}(X)$ has finite injective dimension $\leq n$, and hence has finite flat dimension $\leq n$ as $\Ggd R\leq n$. On the other hand, by \prpcite[7.2]{Gil16} the pair $(^{\perp}{\sf{dw}\widetilde{\sf Inj}},\sf{dw}\widetilde{\sf Inj})$ is an injective cotorsion pair in $\Ch(R)$. Then it follows from Gillespie \corcite[3.3]{Gil17} that $X$ is in $^{\perp}{\sf{dw}\widetilde{\sf Inj}}$, and so is $\sf{\Sigma}X$. Thus one has $\Ext^{1}_{\Ch(R)}({\sf\Sigma}X,I)=0$. This yields that the exact sequence $0 \to I \to \sf{Cone}\alpha \to {\sf\Sigma}X \to 0$ is split. So $\alpha$ is null homotopic; see Enochs, Jenda and Xu \lemcite[3.2]{EJX96}.
\end{prf*}

\begin{cor}\label{main-id}
Let $R$ be a ring with $\GIacgd R$ finite. Then
\begin{equation*}
\underline{\GIac}\simeq\catD{R}/\catK{\sf Inj}\simeq\catD{R}/\catK{\sf Prj}
\end{equation*}
are compactly generated.
\end{cor}
\begin{prf*}
The first equivalence in the statement holds by Theorem \ref{equivalence-dual}, and the second one follows from \thmcite[6.9]{Beligiannis2000} as $\Ggd R<\infty$. Next we prove that $\underline{\GIac}$ is compactly generated. Let ${\sf S}(\sf ACInj)$ denote the homotopy category of all acyclic AC-injective $R$-complexes, and let $\sf{K}_{\sf ac}(\sf Inj)$ (resp., $\sf{K}_{\sf tac}(\sf Inj)$) denote the homotopy category of acyclic (resp., totally acyclic) complexes of injective $R$-modules. Consider the following equivalences:
  \begin{equation*}
    \underline{\GIac}=\underline{\mathsf{GI}}\simeq \sf{K}_{\sf tac}(\sf Inj)=\sf{K}_{\sf ac}(\sf Inj)={\sf S}(\sf ACInj)\:.
  \end{equation*}
Here the first equality holds by Lemma \ref{abclean}. Since $\Ggd R$ is finite by \ref{2.7.1}, all $R$-module have finite Gorenstein injective dimension by Lemma \ref{BMth1.1}. It follows that every acyclic complex of injective $R$-modules has Gorenstein injective cycles and so it is totally acyclic. This yields that the second equality holds. The last equality follows from Lemma \ref{kac-inj}; while the equivalence holds by Krause \prpcite[7.2]{HKr05}. Finally, from \thmcite[5.8 and 4.6]{Gi17} that ${\sf S}(\sf ACInj)$ is compactly generated.
\end{prf*}

Let $R$ be a ring with $\GIacgd R$ finite. Then one has $\underline{\sf GP}\simeq\underline{\sf GI}$ by \thmcite[6.9]{Beligiannis2000} as $\Ggd R<\infty$; see \ref{global}. On the other hand, by Lemma \ref{abclean}, the equality $\GIac={\sf GI}$ holds. So the next result is immediate by Theorem \ref{exa-pd} and Corollary \ref{main-id}.

\begin{cor}\label{main-pd and fd-coherent 2}
Let $R$ be a left coherent ring with $\Ggd R<\infty$. Then
\begin{equation*}
\catD{R}/\catK{\sf Inj}\simeq\catD{R}/\catK{\sf Prj}\simeq\underline{\sf GP}\simeq\underline{\sf GI}
\end{equation*}
are compactly generated.
\end{cor}

We close this section with the following example; it shows that coherent rings of finite global Gorenstein dimension may not be Iwanaga-Gorenstein nor perfect nor Morita\footnote{See \cite{Be01} for the definition of Morita rings. It is known that a ring $R$ is left Morita if and only if $R$ is left Artinian and $\M{R}$ has a finitely generated injective cogenerator.}. Let $R=\prod_{i=1}^{n}R_i$ be a direct product of rings (see \ref{product}). It is easy to see that $R$ is Iwanaga-Gorenstein (resp., left perfect and left Morita) if and only if each $R_i$ is Iwanaga-Gorenstein (resp., left perfect and left Morita).

\begin{exa}\label{exa-2}
Let $R=\mathbb{Z}$ and $S=\left(
         \begin{array}{cc}
           \mathbb{Q} & \mathbb{R} \\
           0 & \mathbb{Q} \\
         \end{array}
       \right)$.
Then $R$ is a commutative Iwanaga-Gorenstein ring that is neither perfect nor Artin; hence $R$ is commutative coherent with $\Ggd R<\infty$.
According to Wang \exacite[3.4]{Wang17} $S$ is a commutative perfect coherent (non-noetherian) ring with $\Ggd S<\infty$. Then the direct product $R\times S$ is a commutative coherent ring with $\Ggd (R\times S)=\sup\{\Ggd R, \Ggd S\}<\infty$, which is neither Iwanaga-Gorenstein nor perfect nor Morita.
\end{exa}

\appendix
\section*{Appendix. Gorenstein AC-flat modules}
\stepcounter{section}
\noindent
In this section we give some properties of Gorenstein AC-flat modules. We notice that all Gorenstein AC-projective $R$-modules are Gorenstein AC-flat. Actually, by \thmcite[A.6]{BGH14}, an $R$-module $M$ is Gorenstein AC-projective if and only if there exists an exact sequence $\cdots\to P_1\to P_0\to P_{-1}\to \cdots$ of projective $R$-modules such that
$M\cong\mathrm{Coker}(P_1\to P_0)$, and it remains exact after applying the functor $A\otimes_R-$ for each absolutely clean $\Rop$-module $A$. The next two results are from Estrada, Iacob and P\'{e}rez \thmcite[2.12]{EIP18} and \exacite[2.17(2)]{EIP18}.

\begin{prp}\label{Gac p and f}
The following conditions are equivalent for an $R$-module $M$.
 \begin{eqc}
 \item $M$ is Gorenstein AC-flat.
\item There is a short exact sequence $0\to K\to L\to M\to0$ of $R$-modules
with $K\in\sf{Flat}$ and $L\in\GPac$, and it remains exact after applying the functor
$\Hom_R(-,C)$ for any (flat) cotorsion $R$-module $C$.
\item $\Ext^1_{R}(M,C)=0$ holds for all cotorsion $R$-modules $C\in (\GPac)^\perp$.
\item There is a short exact sequence $0\to M \to F \to N \to 0$ of $R$-modules
with $F\in\sf{Flat}$ and $N\in\GPac$.
 \end{eqc}
\end{prp}

\begin{prp}\label{Gac p and f1}
The pair $(\GFac, (\GPac)^{\perp}\cap \sf{Cot})$ is a
complete and hereditary cotorsion pair with the equality
$\GFac\cap(\GPac)^{\perp}=\sf{Flat}$.
\end{prp}

Recall that a triple $(\sf{Q},\sf{W},\sf{R})$ of classes of $R$-modules is Hovey triple if $\sf{W}$ is thick
and $(\sf{Q}\cap \sf{W}, \sf{R})$ and $(\sf{Q}, \sf{W}\cap \sf{R})$ are complete cotorsion pairs. If furthermore the above two cotorsion pairs are hereditary then the Hovey triple $(\sf{Q},\sf{W},\sf{R})$ is called hereditary. From Hovey \thmcite[2.2]{Ho02}, an abelian model structure on $\M{R}$ is equivalent to a Hovey triple. This fact is known as ``Hovey correspondence'' in the literature. We hence always denote an abelian model structure $\mathcal{M}$ as a Hovey triple $\mathcal{M}=(\sf{Q},\sf{W}, \sf{R})$.

For an abelian model structure $\mathcal{M}=(\sf{Q},\sf{W}, \sf{R})$, we denote by Ho$(\mathcal{M})$ the homotopy category of $\mathcal{M}$.
By Gillespie \cite[sec 4 and 5]{Gil162}, for any hereditary Hovey triple $\mathcal{M}=(\sf{Q},\sf{W},\sf{R})$, there
is a Frobenius exact category $\sf{Q}\cap \sf{R}$ whose projective-injective objects are precisely those in $\sf{Q}\cap \sf{R}\cap\sf{W}$. Furthermore, the stable category $\underline{\sf{Q}\cap \sf{R}}$ is triangulated equivalent to $\Ho(\mathcal{M})$. This triangulated equivalence is known as the fundamental theorem of model categories in the literature.

\begin{ipg}\label{Frobenius exact category}
By \cite{BGH14}, the triple $\mathcal{M}=(\GPac, \GPac^{\perp}, \M{R})$
is a hereditary Hovey triple. As an immediate consequence of Proposition \ref{Gac p and f1} one gets that the triple $\mathcal{M}'=(\GFac, \GPac^\perp, \sf{Cot})$ is a hereditary Hovey triple, which can also be found in \corcite[4.3]{EIP18}. Thus the category $\GPac$ (resp., $\GFac\cap\sf{Cot}$) is a Frobenius category with projective-injective objects all projective (resp., flat-cotorsion) $R$-modules. By the fundamental theorem of model categories, $\underline{\GPac}$ (resp., $\underline{\GFac\cap \sf{Cot}}$) is triangulated equivalent to ${\rm Ho}(\mathcal{M})$ (resp., ${\rm Ho}(\mathcal{M}')$).
\end{ipg}

It follows from Estrada and Gillespie \lemcite[5.4]{EG19} that if two hereditary Hovey triples $\mathcal{M} =(\sf{Q}, \sf{W}, \sf{R})$ and $\mathcal{M}' = (\sf{Q}', \sf{W}, \sf{R}')$ on $\M{R}$ have the same class $\sf{W}$ of trivial objects and if $\sf{Q} \subseteq \sf{Q'}$ (or equivalently, $\sf{R'} \subseteq \sf{R}$), then there is a triangulated equivalence ${\rm Ho}(\mathcal{M})\simeq{\rm Ho}(\mathcal{M}')$. Applying this fact to the hereditary Hovey triples $\mathcal{M}=(\GPac, \GPac^\perp, \M{R})$ and $\mathcal{M}'=(\GFac, \GPac^\perp, \sf{Cot})$ (in view of \ref{Frobenius exact category}) we get

\begin{cor}\label{triequivalence}
There exits a triangulated equivalence $\underline{\GPac}\simeq\underline{\GFac\cap \sf{Cot}}$.
\end{cor}

Note that the pairs $({\sf Flat},{\sf Cot})$ and $(\GFac, (\GPac)^{\perp}\cap \sf{Cot})$ are complete hereditary cotorsion pairs with $\GFac\cap(\GPac)^{\perp}\cap{\sf Cot}=\sf{Flat}\cap{\sf Cot}$; see Proposition \ref{Gac p and f1}. The following result is immediate by \thmcite[4.5]{DLYZ18}.

\begin{cor}\label{Gac p and f3}
There exists a triangle equivalence
\begin{equation*}
\underline{\GFac\cap \sf{Cot}}\simeq \catD{R}_{\widehat{\GFac}\cap\widehat{\sf{Cot}}}/\catK{\sf FlatCot}.
\end{equation*}
\end{cor}

\section*{Acknowledgments}
\noindent
We thank James Gillespie for conversations and comments on an early draft of this paper, and
we extend our gratitude to the referee for valuable comments that have improved the presentation at several points.

\bibliographystyle{amsplain-nodash}

\begin{thebibliography}{10}

\bibitem{MAsMBr69}
Maurice Auslander and Mark Bridger, \emph{Stable module theory}, Memoirs of the
  American Mathematical Society, No. 94, American Mathematical Society,
  Providence, R.I., 1969. \MR{MR0269685}

\bibitem{BMPS19}
V\'{\i}ctor Becerril, Octavio Mendoza, Marco~A. P\'{e}rez, and Valente
  Santiago, \emph{Frobenius pairs in abelian categories. {C}orrespondences with
  cotorsion pairs, exact model categories, and {A}uslander-{B}uchweitz
  contexts}, J. Homotopy Relat. Struct. \textbf{14} (2019), no.~1, 1--50.
  \MR{MR3913970}

\bibitem{Beligiannis2000}
Apostolos Beligiannis, \emph{The homological theory of contravariantly finite
  subcategories: {A}uslander-{B}uchweitz contexts, {G}orenstein categories and
  (co-)stabilization}, Comm. Algebra \textbf{28} (2000), no.~10, 4547--4596.
  \MR{MR1780017}

\bibitem{Be01}
Apostolos Beligiannis, \emph{Homotopy theory of modules and {G}orenstein rings}, Math. Scand.
  \textbf{89} (2001), no.~1, 5--45. \MR{MR1856980}

\bibitem{BR07}
Apostolos Beligiannis and Idun Reiten, \emph{Homological and homotopical
  aspects of torsion theories}, Mem. Amer. Math. Soc. \textbf{188} (2007),
  no.~883, viii+207. \MR{MR2327478}

\bibitem{DB09}
Driss Bennis, \emph{Rings over which the class of {G}orenstein flat modules is
  closed under extensions}, Comm. Algebra \textbf{37} (2009), no.~3, 855--868.
  \MR{MR2503181}

\bibitem{BM10}
Driss Bennis and Najib Mahdou, \emph{Global {G}orenstein dimensions}, Proc.
  Amer. Math. Soc. \textbf{138} (2010), no.~2, 461--465. \MR{MR2557164}

\bibitem{BEI18}
Daniel Bravo, Sergio Estrada, and Alina Iacob, \emph{{\rm FP}$_n$-injective, {\rm FP}$_n$-flat covers and preenvelopes, and Gorenstein AC-flat covers}, Algebra Colloq. \textbf{25} (2018), no.~2, 319--334. \MR{MR3805326}

\bibitem{BGH14}
Daniel Bravo, James Gillespie, and Mark Hovey, \emph{{The stable module
  category of a general ring}}, {preprint,} \textbf{\arxiv[RA]{1405.5768}}.

\bibitem{ROB86}
Ragnar-Olaf Buchweitz, \emph{Maximal {C}ohen--{M}acaulay modules and
  {T}ate-cohomology over {G}orenstein rings}, University of Hannover, 1986,
  available at \mbox{\sffamily http://hdl.handle.net/1807/16682}.

\bibitem{Chen11}
Xiao-Wu Chen, \emph{Relative singularity categories and {G}orenstein-projective
  modules}, Math. Nachr. \textbf{284} (2011), no.~2-3, 199--212. \MR{MR2790881}

\bibitem{CET19}
Lars~Winther Christensen, Sergio Estrada, and Peder Thompson, \emph{{Homotopy
  categories of totally acyclic complexes with applications to the
  flat--cotorsion theory}}, {Contemp. Math, to appear,}
  \textbf{\arxiv[RA]{1812.04402v2}}.

\bibitem{DLYZ18}
Zhenxing Di, Zhongkui Liu, Xiaoyan Yang, and Xiaoxiang Zhang,
  \emph{Triangulated equivalence between a homotopy category and a triangulated
  quotient category}, J. Algebra \textbf{506} (2018), 297--321. \MR{MR3800079}

\bibitem{IEm12}
Ioannis Emmanouil, \emph{On the finiteness of {G}orenstein homological
  dimensions}, J. Algebra \textbf{372} (2012), 376--396. \MR{MR2990016}

\bibitem{EEnOJn95b}
Edgar~E. Enochs and Overtoun M.~G. Jenda, \emph{Gorenstein injective and
  projective modules}, Math. Z. \textbf{220} (1995), no.~4, 611--633.
  \MR{MR1363858}

\bibitem{EJT-93}
Edgar~E. Enochs, Overtoun M.~G. Jenda, and Blas Torrecillas, \emph{Gorenstein
  flat modules}, Nanjing Daxue Xuebao Shuxue Bannian Kan \textbf{10} (1993),
  no.~1, 1--9. \MR{MR1248299}

\bibitem{EJX96}
Edger~E. Enochs, Overtoun M.~G. Jenda, and Jinzhong Xu, \emph{Orthogonality in
  the category of complexes}, Math. J. Okayama Univ. \textbf{38} (1996),
  25--46. \MR{MR1644453}

\bibitem{EG19}
Sergio Estrada and James Gillespie, \emph{The projective stable category of a
  coherent scheme}, Proc. Roy. Soc. Edinburgh Sect. A \textbf{149} (2019),
  no.~1, 15--43. \MR{MR3922806}

\bibitem{EIP18}
Sergio Estrada, Alina Iacob, and Marco A. P\'{e}rez, \emph{Model structures and relative Gorenstein flat modules and chain complexes}, Contemporary Math., {to appear,} \textbf{\arxiv[RT]{1709.00658v2}}.

\bibitem{Gi04}
James Gillespie, \emph{The flat model structure on {${\rm Ch}(R)$}}, Trans.
  Amer. Math. Soc. \textbf{356} (2004), no.~8, 3369--3390. \MR{MR2052954}

\bibitem{Gil16}
James Gillespie, \emph{Gorenstein complexes and recollements from cotorsion pairs},
  Adv. Math. \textbf{291} (2016), 859--911. \MR{MR3459032}

\bibitem{Gil162}
James Gillespie, \emph{Hereditary abelian model categories}, Bull. Lond. Math. Soc.
  \textbf{48} (2016), no.~6, 895--922. \MR{MR3608936}

\bibitem{Gil17}
James Gillespie, \emph{On {D}ing injective, {D}ing projective and {D}ing flat modules
  and complexes}, Rocky Mountain J. Math. \textbf{47} (2017), no.~8,
  2641--2673. \MR{MR3760311}

\bibitem{Gi17}
James Gillespie, \emph{On the homotopy category of {AC}-injective complexes}, Front.
  Math. China \textbf{12} (2017), no.~1, 97--115. \MR{MR3579262}

\bibitem{Gil19}
James Gillespie, \emph{A{C}-{G}orenstein rings and their stable module categories}, J.
  Aust. Math. Soc. \textbf{107} (2019), no.~2, 181--198. \MR{MR4001567}

\bibitem{Ho02}
Mark Hovey, \emph{Cotorsion pairs, model category structures, and
  representation theory}, Math. Z. \textbf{241} (2002), no.~3, 553--592.
  \MR{MR1938704}

\bibitem{KK88}
Ellen Kirkman and James Kuzmanovich, \emph{On the global dimension of fibre
  products}, Pacific J. Math. \textbf{134} (1988), no.~1, 121--132. \MR{MR953503}

\bibitem{Krause00}
Henning Krause, \emph{Smashing subcategories and the telescope conjecture---an
  algebraic approach}, Invent. Math. \textbf{139} (2000), no.~1, 99--133.
  \MR{MR1728877}

\bibitem{HKr05}
Henning Krause, \emph{The stable derived category of a {N}oetherian scheme}, Compos.
  Math. \textbf{141} (2005), no.~5, 1128--1162. \MR{MR2157133}

\bibitem{MD06}
Lixin Mao and Nanqing Ding, \emph{The cotorsion dimension of modules and
  rings}, Abelian groups, rings, modules, and homological algebra, Lect. Notes
  Pure Appl. Math., vol. 249, Chapman \& Hall/CRC, Boca Raton, FL, 2006,
  pp.~217--233. \MR{MR2229114}

\bibitem{Orlov04}
Dmitri Orlov, \emph{Triangulated categories of singularities and {D}-branes in
  {L}andau-{G}inzburg models}, Tr. Mat. Inst. Steklova \textbf{246} (2004),
  no.~Algebr. Geom. Metody, Svyazi i Prilozh., 240--262. \MR{MR2101296}

\bibitem{SS18}
Jan \v{S}aroch and Jan \v{S}\v{t}ov\'{\i}\v{c}ek, \emph{{Singular compactness
  and definability for $\Sigma$-cotorsion and Gorenstein modules}}, Selecta Math. \textbf{26} (2020), article number: 23.

\bibitem{Wang17}
Junpeng Wang, \emph{Ding projective dimension of {G}orenstein flat modules},
  Bull. Korean Math. Soc. \textbf{54} (2017), no.~6, 1935--1950. \MR{MR3733774}

\bibitem{YD15}
Xiaoyan Yang and Nanqing Ding, \emph{On a question of {G}illespie}, Forum Math.
  \textbf{27} (2015), no.~6, 3205--3231. \MR{MR3420339}

\end{thebibliography}

\def\cprime{$'$}
  \providecommand{\arxiv}[2][AC]{\mbox{\href{http://arxiv.org/abs/#2}{\sf
  arXiv:#2 [math.#1]}}}
  \providecommand{\oldarxiv}[2][AC]{\mbox{\href{http://arxiv.org/abs/math/#2}{\sf
  arXiv:math/#2
  [math.#1]}}}\providecommand{\MR}[1]{\mbox{\href{http://www.ams.org/mathscinet-getitem?mr=#1}{#1}}}
  \renewcommand{\MR}[1]{\mbox{\href{http://www.ams.org/mathscinet-getitem?mr=#1}{#1}}}
\providecommand{\bysame}{\leavevmode\hbox to3em{\hrulefill}\thinspace}
\providecommand{\MR}{\relax\ifhmode\unskip\space\fi MR }
\providecommand{\MRhref}[2]{%
  \href{http://www.ams.org/mathscinet-getitem?mr=#1}{#2}
}
\providecommand{\href}[2]{#2}

\end{document}